\documentclass[twoside,11pt,preprint]{article}

\usepackage{jmlr2e}

\usepackage{microtype}
\usepackage{booktabs,bbm,amsmath,xspace}
\usepackage{xcolor,bm,amsmath,graphicx,comment,bigints}
\usepackage{algorithm,algorithmic,caption,subcaption}

\newcommand\cvh{$\mathcal{H}$\xspace}
\newcommand\dtr{$\mathcal{D}$\xspace}

\ShortHeadings{}{}
\firstpageno{1}

\begin{document}

\title{A Sketching Method for Finding the Closest Point on a Convex Hull}

\author{\name Roozbeh Yousefzadeh \email roozbeh.yousefzadeh@yale.edu \\
       \addr Yale Center for Medical Informatics and VA Connecticut Healthcare System
       }


\maketitle

\begin{abstract}
We develop a sketching algorithm to find the point on the convex hull of a dataset, closest to a query point outside it. Studying the convex hull of datasets can provide useful information about their geometric structure and their distribution. Many machine learning datasets have large number of samples with large number of features, but exact algorithms in computational geometry are usually not designed for such setting. Alternatively, the problem can be formulated as a linear least-squares problem with linear constraints. However, solving the problem using standard optimization algorithms can be very expensive for large datasets. Our algorithm uses a sketching procedure to exploit the structure of the data and unburden the optimization process from irrelevant points. This involves breaking the data into pieces and gradually putting the pieces back together, while improving the optimal solution using a gradient projection method that can rapidly change its active set of constraints. Our method eventually leads to the optimal solution of our convex problem faster than off-the-shelf algorithms.
\end{abstract}


\section{Introduction} \label{sec:intro}

Studying the convex hull of datasets can provide useful information about their geometric structure and their distribution. Such information may then help with downstream tasks such as out-of-distribution detection. Studying the relationship of testing sets with respect to the convex hull of training sets may provide insights about the generalization of models, too \citep{yousefzadeh2021convexhull}.



Standard convex hull algorithms in computational geometry are highly tailored and efficient, but most of them are not practical in high-dimensional domains. There are approximation algorithms as well with useful guarantees \citep{blum2019sparse}, but sometimes we may need to seek exact solutions.

Here, we provide a fast algorithm to understand the geometric relation of a query point, $q$, in high-dimensional space with respect to the convex hull of a dataset, \cvh, possibly with large number of datapoints. Specifically, we find the point on the body of the convex hull, closest to the query point outside it. The vector connecting the query point to the convex hull reveals useful information about the query point and also about the contents of the dataset. The size of the vector in comparison with the size of \cvh tells us how far $q$ is from the distribution of \cvh. Beyond the size, that vector reveals the direction that can bring the $q$ to the convex hull. This also tells us something about the dataset. The closest point on the convex hull is a convex combination of certain number of points in the dataset. Investigating those points and their geometric relationship w.r.t. $q$ can reveal further information about the data.

\section{Formulation} \label{sec:formulation}


Let's consider that dataset \dtr is formed as a matrix, with $n$ rows corresponding to the samples, and $d$ columns corresponding to the features, i.e., dimensions. \cvh is the convex hull of all the samples in \dtr. Our query point, $q$, sits outside the \cvh, and we seek to find $x^{\mathcal{H}}$, the point in \cvh that is closest to $q$.

To ensure that $x$ belongs to \cvh, we can define
\begin{equation} \label{eq:alpha}
    x^{\mathcal{H}} = \alpha \mathcal{D},
\end{equation}

where $\alpha$ is a row vector of size $n$. If all elements of $\alpha$ are bounded between 0 and 1, and their summation also equals 1, then by definition, $x$ belongs to \cvh. Given equation~\eqref{eq:alpha}, we can change our optimization variable to $\alpha$.

Our objective function is:

\begin{equation} \label{eq:obj}
    \min_{\alpha} f(\alpha) = \| q - \alpha \mathcal{D}\|_2^2,
\end{equation}

while our constraints ensure that $x$ belongs to \cvh.

\begin{equation} \label{eq:const1}
    \alpha \mathbbm{1}_{n,1} = 1,
\end{equation}
\begin{equation} \label{eq:const2}
    0 \leq {\alpha}.
\end{equation}

This is a constrained least squares problem which can be solved using standard algorithms in numerical optimization literature \citep{nocedal2006numerical}. For any query point, we first compute the optimal $\alpha$ using equations \eqref{eq:obj}-\eqref{eq:const2}. We then compute the corresponding $x^{\mathcal{H}}$ using the optimal $\alpha$ and equation~\eqref{eq:alpha}.

However, we note that \dtr can be quite large and most likely, a large portion of points in \cvh would not have any effect on the optimal solution. This is easy to envision in 2D. Consider, for example, a large set of points forming a square convex hull, and a single point $q$, outside it. Now assume that one edge of the square is closer to $q$, compared to the other three edges. Then, for our optimization problem, we can only use the points on that closer edge. Whether we include the other points or not, the solution to our optimization problem will remain the same. Hence, we can exploit the geometric structure of \dtr with respect to $q$, and solve our optimization problem faster.

\section{Our Algorithm}


Here, we develop an algorithm that is suitable for datasets with large number of samples in high-dimensional space.

\subsection{Preliminaries}

We expect the optimal solution to be highly sparse. For example for a dataset with $n=60,000$, the optimal solution may have only $100$, or even $30$ non-zero elements. 
Therefore, if we identify the subspace that contains our optimal solution, finding the optimal solution in that subspace will be easy. This observation motivates us to start by solving our optimization problem for a subset of dataset denoted by $\mathcal{D}^\prime$, and then gradually enlarge that subset until it includes the entire~$\mathcal{D}$.

Let's denote the convex hull of \dtr by \cvh. Let's also assume that the starting point for our optimization problem is some point $x^0$ inside the \cvh. This means that $x^0$ satisfies the constraints~\eqref{eq:const1}-\eqref{eq:const2}. As the optimization algorithm makes progress, the solution gradually moves towards the $q$ until it reaches a point where it cannot move any closer to $q$ without exiting the \cvh. That point, denoted by $x^{\mathcal{H}}$, is the optimal solution we seek to find.

To develop our algorithm, let's consider a subset of $n'$ samples from \dtr, and a specific $q$. The optimal solution of \eqref{eq:obj}, subject to \eqref{eq:const1}-\eqref{eq:const2}, for the subset, would be $\alpha^\prime$, leading to closest point $x^\prime$. If some element $i$ of $\alpha^\prime$ is zero, it means that the sample $i$ has no effect on the optimal solution $\alpha^\prime$, and by extension, it has no effect on the $x^\prime$. Hence, if we discard sample $i$ from the beginning, $x^\prime$ will not change. However, if we add an additional point to the subset, so that the size of the set increases to $n^\prime+1$, sample $i$ might become part of the optimal solution. Hence, the subspace that contains $x^\prime$ may include dimensions that do not belong to the subspace of $x^{\mathcal{H}}$. Moreover, there may be dimensions included in the subspace of $x^{\mathcal{H}}$ that are not included in $x^\prime$.

Therefore, as we increase the size of $\mathcal{D}^\prime$, we have to consider all the included samples because any of them may become part of the optimal solution. At each iteration, the solution is trying to get closer to $q$ by moving inside a convex body. Constraint~\eqref{eq:const1} is equality and will be satisfied as long as the solution stays within the \cvh. We expect the lower bound of constraint~\eqref{eq:const2} to be binding for many of the samples, but we do not expect it to be upper bounded unless the optimal solution coincides with a point in the dataset, in which case, constraint~\eqref{eq:const2} will become binding for all the samples (upper bounded for that point and lower bounded for all other points in $\mathcal{D}$).



\subsection{Gradient Projection Method}

To move from $x^0$ to the $x^{\mathcal{H}}$, we can use the Gradient Projection Method described by \citet[Section 16.7]{nocedal2006numerical} which we briefly review in the following.

At each iteration of the Gradient Projection Method, we first compute the Cauchy point (the feasible minimizer along the direction of derivatives), and then perform a subspace minimization for the samples that are not binding in the lower bound of constraint~\eqref{eq:const2}.

For computing the Cauchy point, we use the same procedure described by \citet[Section 16.7]{nocedal2006numerical}. To make sure that constraint~\eqref{eq:const1} is satisfied and the Cauchy point does not exit the convex hull region, we normalize the gradient direction so that the sum of its values equals zero.

The subspace minimization that follows does not need to lead to an exact solution of the sub-problem, as it could make the process unnecessarily expensive. As long as it improves the solution at hand, we can move to the next iteration and compute a new Cauchy point.

\subsection{Dual form}

We note that our constraint~\eqref{eq:const1} is an equality constraint. The dual form of our problem only has non-negativity constraint on the Lagrange multipliers. The gradient projection method can be applied to the dual form of our problem, and in certain situations, solving the dual form may be easier. 

In practice, for our datasets of interest, we observe that Lagrange multipliers are mostly non-sparse and that makes the dual form more expensive to solve. Here, we present the dual form of the problem because it may be useful in certain settings when the primal solution is non-sparse.

The Lagrangian for our optimization problem is
\begin{equation} \label{eq:lagrangian}
    \mathcal{L}(\alpha,\lambda) = \frac{1}{2} \| {q} - {\alpha} {\mathcal{D}}\|_2^2 - ({\alpha \mathbbm{1}}_{n,1} - 1) \lambda^1 - \alpha \lambda^2,
\end{equation}
where $\lambda$'s are Lagrange multipliers: $\lambda^1$ is a scalar while $\lambda^2$ is a column vectors with $n$ elements.

The Lagrange dual objective is
\begin{equation} \label{eq:dual_obj}
    g(\lambda) = \inf_{\alpha} \; \mathcal{L}(\alpha,\lambda).
\end{equation}

Because $\mathcal{L}(.,\lambda)$ is a strictly convex quadratic function, the infimum is achieved when $\nabla_{\alpha} \mathcal{L}(\alpha,\lambda) = 0$, i.e.,

\begin{equation} \label{eq:dual_opt}
    - \mathcal{D} q^T + \mathcal{D} \mathcal{D}^T \alpha^T  - \lambda^1 \mathbbm{1}_{n,1} - {\lambda^2}  = 0.
\end{equation}

If we compute the Singular Value Decomposition of $\mathcal{D} = U \Sigma V^T$ and plug it into \eqref{eq:dual_opt}, we obtain

\begin{equation*} \label{eq:dual_opt_svd}
    - U \Sigma V^T q^T + {U} \Sigma^2 {U}^T \alpha^T  - \lambda^1 \mathbbm{1}_{n,1} - {\lambda^2}  = 0,
\end{equation*}

which leads to

\begin{equation} \label{eq:dual_opt_alpha}
    \alpha^d = 
    q V \Sigma^{-1} U^T + (\lambda^1 \mathbbm{1}_{1,n} + {\lambda^2}^T)U \Sigma^{-2} U^T.
\end{equation}

Then, the dual form of our problem is
\begin{equation} \label{eq:dual}
    \max_{\lambda \in \mathbbm{R}} \; g(\lambda) = \mathcal{L}(\alpha^d,\lambda)  \quad \text{subject to:} \quad \lambda^1, \lambda^2 \geq 0.
\end{equation}

At each iteration of the Gradient Projection Method, we first compute $\nabla_\lambda g(\lambda)$, then find the Cauchy point in the direction of gradient, and perform an inexact minimization in the subspace of Lagrange multipliers that are not in the active set. This process repeats until the KKT conditions are satisfied, at which point, the corresponding $\alpha^d$ will be our optimal solution.

For subspace minimization, one can use the alternating direction method of multipliers (ADMM) on $\lambda^1$ and $\lambda^2$.



\subsection{Sketching Algorithm}

Now, we have all the pieces to formalize our sketching method in Algorithm~\ref{alg:cvx}.

Our sketching method divides the dataset, \dtr, into $\eta$ pieces: $\mathcal{D}_1,\mathcal{D}_2, \dots, \mathcal{D}_\eta$. It initiates $\Phi$ by adding only the first piece of data, $\mathcal{D}_1$. It then solves equations~\eqref{eq:obj}-\eqref{eq:const2} for $\mathcal{D}_1$ using the gradient projection method described before. When it finds the optimal solution, $\alpha^*$, it proceeds with appending $\mathcal{D}_2$ to the $\Phi$, and solves the problem again using the optimal solution from previous step. We expect the solution from previous step to be sparse implying that its active set is relatively large. We also expect the new Cauchy point to not be much different in the subspace of $\mathcal{D}_1$ and we expect many of the active sets to remain active in that subspace. 

This is the key benefit of sketching, because computing the Cauchy point is not expensive, and we have already excluded a considerable portion of $\mathcal{D}_1$ from the expensive part of computations. And this benefit repeats at the next sketching step because a considerable portion of $\mathcal{D}_2$ will be included in the active set and thereby excluded from the following subspace minimization.

In other words, the subset of points that we pick at the beginning for solving the problem provide a relatively good sketch of the entire convex hull w.r.t the $q$. And we gradually improve our sketch until we obtain the full picture of the \cvh. The sketching algorithm ends when we have included all pieces of $\mathcal{D}$ in $\Phi$ and we obtain the optimal solution for $\mathcal{D}$.

Algorithm~\ref{alg:cvx} formalizes this entire process.

\begin{algorithm}[H]
\caption{Finding the point on the convex hull of a dataset, closest to a query point outside it}
\label{alg:cvx}
\textbf{Inputs}: Dataset \dtr, query point $q$, number of partitions $\eta$\\
\textbf{Outputs}: $x^*$: the point on the convex hull of \dtr, closest to $q$
\begin{algorithmic}[1] 
\STATE Sort the rows in \dtr based on their closeness to $q$
\STATE Partition \dtr into $\eta$ subsets, call each subset $\mathcal{D}_i$ ($\mathcal{D}_1$ contains the points closest to $q$)
\STATE Initialize an empty matrix $\Phi$
\STATE Initialize the optimal solution $\alpha^*$ as an empty vector
\FOR{$i=1$ to $\eta$}
    \STATE Append $\mathcal{D}_i$ to $\Phi$
    \STATE Append zeros to vector $\alpha^*$ for each sample in $\mathcal{D}_i$
    \WHILE{KKT conditions are not satisfied for objective function \eqref{eq:obj}, subject to  constraints \eqref{eq:const1}-\eqref{eq:const2} on $\Phi$}
        \STATE Compute $\nabla_\alpha f(\alpha)$
        \STATE Compute the Cauchy point in the direction of $\nabla_\alpha f(\alpha)$
        \STATE Approximately solve equation \eqref{eq:obj}, subject to  \eqref{eq:const1}-\eqref{eq:const2} in the subspace of inactive constraints
    \ENDWHILE
\ENDFOR
\STATE $x^* = \alpha^* \Phi$ 
\STATE \textbf{return} $x^*$
\end{algorithmic}
\end{algorithm}

\acks{R.Y. thanks Daniel Robinson for his nonlinear optimization course. R.Y. was supported by a fellowship from the Department of Veterans Affairs. The views expressed in this manuscript are those of the author and do not necessarily reflect the position or policy of the Department of Veterans Affairs or the United States government.}



\vskip 0.2in
\bibliography{refs}

\end{document}